\documentclass{article}

\usepackage[affil-it]{authblk}
\usepackage{amsmath, hyperref,cleveref, caption,floatrow, enumerate}

\DeclareNewFloatType{exfloat}%
	{name=Example}

\DeclareNewFloatType{nonexfloat}%
	{name=Non-Example}

\newcommand\mycaption{}
\newcommand\mylabel{}

\newenvironment{example}[2]
	{
	\renewcommand{\mycaption}{#1}%
	  \renewcommand{\mylabel}{#2}%
	\begin{exfloat}[h]
	\centering
	\begin{tabular}{ p{30em} }
	\hline
	\hline\\
	}
	{
	\\\\\hline
	\hline
	\end{tabular}
	\caption{\mycaption}\label{ex::\mylabel}%
	\end{exfloat}
	}

\newenvironment{nonexample}[2]
	{
	\renewcommand{\mycaption}{#1}%
	  \renewcommand{\mylabel}{#2}%
	\begin{nonexfloat}[h]
	\centering
	\begin{tabular}{ p{30em} }
	\hline
	\hline\\
	}
	{
	\\\\\hline
	\hline
	\end{tabular}
	\caption{\mycaption}\label{nex::\mylabel}%
	\end{nonexfloat}
	}

\title{Improv practices in Mathematics Active Teaching}
\author{Pedro Morales-Almazan\\ email \href{mailto:pmorale5@ucsc.edu}{pmorale5@ucsc.edu} }

\affil{Mathematics Department,\\ University of California, Santa Cruz}

\date{\today}

\begin{document}

\maketitle

\begin{abstract}
We explore the parallelism between theater Improv and teaching in an Active Learning environment. We present the notions of Active Teaching as a natural complement of Active Learning, and how unexpected situations give rise to valuable Teaching Moments. These Teaching Moments can be strategically utilized following the rules of Improv. We present some example of this in the Mathematics classroom, as well as the implementation of an Improv Seminar in the Mathematics Department at the University of California, Santa Cruz.
\end{abstract}

\section{Introduction}

Active learning has become an alternative to traditional lecture in making the student a more active part of the learning experience. Specially in mathematics, active learning plays a crucial role in motivating the students. 

Generally speaking, the term Active Learning deals with the student becoming an active component of class. As a natural consequence, class might deviate from scripts to adjust to the student interaction and participation. This would suggest a more free form in class were the instructor should embrace certain unexpected situations. By increasing the interaction with and within students, we also increase the unexpected. This is the reason why Active Learning also requires a degree of Active Teaching. By this we mean the quality of making the teacher an active component of the class together with the student. 

Knowing how to handle the unexpected is an essential skill in teaching. Unexpected situations generally carry naturally occurring Learning Moments. When the teacher becomes active, she can identify these moments and incorporate them into the class. Using spontaneous learning moments is a great way to motivate content in class without forcing the material. By responding to the present action of the student, the teacher can improve the communication, and overall, the learning experience. 

The active teacher is required to handle the unexpected, or to improvise, in many situations. Sometimes the term \emph{improvisation} could have a negative connotation, denoting something of bad quality or a temporal fix without a preparation. However, improvisation in music,  theater, and dance has become an art form in itself. Despite not being previously scripted, improvisation requires extensive training and preparation. It also enhances the creative and emotional parts in the performances.

Improvisation in the arts can be taught and mastered in a systematical way. The same guiding principles used to explore creativity and emotion in improvisation can also be used in teaching. These come as a natural complement to the spontaneity obtained in the student's active learning. 

In this paper, we explore improvisation as a systematic art form and how its guiding principles can be applied to the instructor mathematics instructor.  Section Two of the paper introduces improvisation in the arts and focuses on theater improvisation, or Improv for short. Then it explores the connections with teaching and how improvisation arises naturally in an active learning environment. After this, a more systematic approach is made by analyzing the four pillars of Improv: building a base reality, acknowledging the other, creating an emotional connection, and being present. These four pillars are explored in the Improv context and then formulated in the teaching environment. Following the four pillars, we describe the rules of Improv, which become a practical way to incorporate them. We present five rules used in Improv and how these can apply to the teaching context. We also present real class examples in mathematics were the rules were applied and how this made a difference on the lesson plan and the students' learning. 

On Section tTree we describe an implementation of an Improv seminar done in the Mathematics Department at the University of California, Santa Cruz. This seminar was held once a week for eight weeks, mainly with Teaching Assistants and Tutors from the Mathematics Department. The approach for this seminar was to use Improv games that explored the rules of Improv together with a discussion on the pedagogical implications of the pillar of Improv. We also present some of the games used in each of the sessions.

Finally, Section Four present the conclusions and closing remarks of the paper. 

\section{Improv}

Improvisation appears in several forms of art. The most common ones are music, dance, and theater. Usually in music, we experience pieces that are carefully rehearsed either in a recorded version or while listening to live music. Also in dance, choreographic presentations like ballet are very well structured and carefully planned. Theater presents these characteristics as well. In all of these performing arts, structure, planning, and rehearsals are fundamental. Artists strive to improve their performances. Top performing artists spend a significant amount of time practicing and rehearsing. Planning and scripting traditionally are an important part of these performing arts. Composers, choreographers, and writers play a crucial role. In the end, it is the fusion of these and the performers that bring presentation to life.

On the other hand, there are improvised versions of these performing arts. Music jams or improvised sessions are common in a more informal setting. Improvised dancing is very common in social events. Improvised theater, or Improv, is very popular in comedy. Improvisation is also often used as a way to warm up for performances. Artists use improvisation to practice performing techniques, to do team-building, and to get into the performance atmosphere.

Improvisation also has powerful educational value. It is used to stimulate creativity, develop flexibility, build team work, and developing self-awareness. The value of these skills is not only restricted to the arts, but to almost any field. Recently, the usage of Improv techniques in general professional settings is becoming more and more common. Companies use Improv workshops to develop team building, stimulate creativity in R\&D divisions, and improve communication flow across departments. 

This appears to the increasing need of develop soft skills in professionals, specially in technical fields such as STEM. These Improv workshops aim to complement professional development from a social perspective, incorporating emotional, creative, social, and communication components to the professional setting.
 
Most of these workshops and seminars use ideas and exercises from improvisational theater. The main idea of improvisational theater, or Improv for short, is to develop a play in real time without a script. Often, cues from the audience are required as the initial inspiration for the play, letting the spontaneity of the actions of the performers dictate the next moves. Like in music and dance, being able to successfully improvise a play requires much preparation and technical skills. Performers require many years of preparation and practice to be able to improvise a scene. 

Conceptually, Improv could be perceived as unexpected and not organized. However, guided or disciplined improvisation enables to take advantage of the skills developed by improvisation in a systematic way. Within the free form of improvisation, it is still possible to apply certain higher level guidelines that keep the unexpected under control and use it to our advantage. The apparent lack of prescribed structure on the performance level is compensated by a set of guidelines in a higher level of abstraction. These are often called \emph{rules of Improv}, which enclose strategies to, despite not having a script, let performers weave a unified and consequential story line.

The parallelism between teaching and artistic performance has been extensively explored \cite{sawyer2004creative, sawyer2011structure, berk2009whose}. We can say that there are two main classifications within this. 

The first is to think about teaching as a scripted play. The instructor designs her course and her lesson plans. As a scripted play, the instructor would carefully think about the right definitions, examples, and activities to include in her lecture. She would script the content seen in class, and to the best of her experience, will foresee the places where students could have trouble. Also, the more experience instructor would be able to identify the type of questions that students might come up with. Hence, during lecture, she can focus on the delivery, similarly to how an actress focuses on her performance on a play, and not on the script. This includes paying attention to her voice, visibility of materials in the board or document camera, her class layout, basically her stage presence.

The second way is to think about it as improvisation. This is a more flexible approach, where the instructor is more attentive to respond to the particular needs and reactions of students. In this case, the instructor pays more attention to the interactions made by students and uses them to her advantage in guiding her agenda for the class. This, as with artistic improvisation, also requires skill and practice to obtain a successful lesson.

Using Improv as a metaphor for teaching offers a way to pursue a more creative teaching that aligns with the emotional state of students \cite{sawyer2004creative, sawyer2011structure, borko1989cognition}. This becomes a very natural approach for addressing discussions and interactions in class. It is often said that \emph{life is improvised}, and as part of life, this applies for everyday conversations and discussions. Thinking about discussions as interaction between two or more people, we can think of the importance that the \emph{rules of Improv} can have in making discussions more effective. In the context of teaching, we face discussions in different ranges. From Socratic seminars with a group of student to small interactions generated by students' questions.  These interactions become very common when following a teaching approach using active-learning, constructivism, or problem-based learning \cite{sawyer2004creative}.

The discussions generated with students constitute a rich source for \emph{teaching moments}. In mathematics, we often fall in the unfortunate habit of teaching answers to questions students never had before. Sometimes it becomes very difficult to motivate concepts in mathematics courses, making it more dry for students to grasp the intuition behind the material. One of the main goals behind \emph{active learning} is to promote and take advantage of these teaching moments. This is often achieved by encouraging interaction among students and with the instructor. These interactions generate discussions that could escape the instructor's original plan and scope, but create unscripted teaching moments. These are important, as students are usually more motivated and receptive to learning when the teaching is not forced, but naturally motivated. 

This \emph{active teaching} requires then the instructor to be receptive of the student emotional an intellectual state. By deviating from her script but reacting to the student, the instructor is being more effective taking advantage of the teaching moments that naturally arise within the student interactions. Improv is often associated with comedy, and often the comic parts of an improvised scene are not forced. Rather, these appear naturally as a consequence of the emotional and intellectual connection among the performers on stage. Similarly, we can say that learning cannot be forced, but promoted. Teaching moments naturally arise when there is a deeper emotional and intellectual connection between students and instructor. Then, the key role of the instructor is to promote, identify, and use these teaching moments. The main idea is to complement active learning with active teaching in order to produce teaching moments.

\subsection{Pillars of Improv}
\label{sec:impro}

Despite the unscripted nature of Improv, it can be given some structure and \emph{rules} to guide it. These commonly called \emph{rules} can be thought of more as techniques or guidelines that can help the improvisers to build a coherent scene. In Improv, we can say that good Improv is obtained when there is a intellectual and emotional connection between the performers. These rules then emerged heuristically as a way to identify patterns that could contribute to this end \cite{napier2015improvise}. 

There is no clear consensus on all the rules of Improv, but there are some key aspects that are usually included. There are four pillars to take into account when improvising:
\begin{enumerate}
\item Build a base reality
\item Acknowledge the other
\item Create an emotional connection
\item Be present
\end{enumerate}

These four pillars are the basis for the rules of Improv in general. These are goals that every performer should keep in mind when creating a scene and that can also be applied to instructors when teaching. 

\subsubsection{Build a base reality}

This pillar refers to \emph{agreement}. Some people call this the \emph{who, what, where}. Who are the characters, what are they doing, and where are they. This corresponds to establishing the shared knowledge about the scene by framing a location, identifying the characters, and specifying a goal or purpose. Before getting into a plot, it is important for the performers to have a context to build a story on. This also constitutes the ground rules for the reality on the scene. If the framing of the scene happens in a beach setting, the performers and the audience can expect talking about water, sun, sand, etc. This also helps the performers to develop an organic plot to the scene. By framing the environment, the key elements of the story can naturally unfold for the performers, and thus seem smooth for the audience. By having the context of the scene in the beach, a plot line could be ``finding a place where to put my umbrella.'' Then characters could play that one of them is trying to put their umbrella but everywhere they try there is someone else already lying down.The context gives a hint of possible plot lines. 

Similarly, when teaching we should frame our lessons with a base reality. This encloses not only the physical environment, but also the intellectual and emotional reality. That is, framing physically where is the lesson taking place, who is the audience, what kind of resources are at hand, what is the background of the students and the instructor. It also refers to framing the ground rules for class. What is expected from students during the lesson, how are interactions managed, what is permitted in class, what resources are students allowed to use, etc. Another important component is the intellectual framing. What situations are we modeling, what kind of questions are we trying to answer, what is the goal of the theory, what can the topics be used for, etc. Finally, framing the emotionally reality entails to recognize the relevance of the topics for the instructor and students on a more personal level. How can the topics be used in their majors, in their own career, how the topics have been used in prior situations by the instructor or the students, why would the instructor or students care about the topics, etc. 

Building a base reality means to reach a common ground of rules, goals, and motivations that can frame the lessons in a more fluid way. Setting the who, what, where becomes a very important component for the lesson.  

\subsubsection{Acknowledge the other}

An important part of improvisation is the \emph{team work} pillar. One of the most challenging parts of improvisation is that it is done with others. It is not a solo performance, but rather, is a team job. Therefore it is important to avoid being self-center and to embrace the other. To acknowledge the other performers means to recognize their needs and wants, and not to deny them. Denial is one of the most dangerous things in Improv. 

When denying an offer from another performer on stage, the audience can perceive a lack of cohesion or harmony in the scene. This can result in forced plots or cyclic arguments. To acknowledge the other performer mean to take their intellectual and emotional states and recognize them as true and valid. For example, if a performer says ``it is very cold, please close that window,'' and the other performer says `` it is not cold, and there are no windows here!'' then the scene falls apart. The audience will see two independent characters on stage that don't have any base reality. This also might discredit the characters as either been lying or being delusional. It is also very hard to continue a plot from there, since it seems that the two characters are on different realities. 

On the other hand, if the second character would reply ``I don't want to close the window because I'm lazy, and I have a sweater, so I'm not cold either'' even thought the character doesn't second the intention of the original statement, it is not denying the reality. It acknowledges the existence of the window and the fact that the first character is cold. The second character is choosing to adopt an indifferent reaction, which is still compatible with the reality presented by the first character. Now the plot can naturally flow to a conflict between the two. The first character can reply back with something like `` you never care of my feelings!'' and the scene can continue in a natural way.

For the instructor, this means to recognize where the student is with respect to the content. A very common situation is when students say ``this is hard'' or `'`this is difficult.'' As instructors we often fall in the unconscious denial of saying ``no, this is easy!'' This denial is saying that the student is not capable. That somehow is missing something \emph{easy}. Like in the performance case, rejecting the students intellectual and emotional place can lead to a learning lock. On the other hand, by recognizing the student's reaction, the instructor is in a good place to pick it up and guide the flow of the conversation, promoting the student's learning. By acknowledging the perception of difficulty for the student, the instructor can connect with the her and work their way through the concept. 

The instructor may reply ``it is challenging, but we will work through it. Why do you say it is difficult?'' Notice that this response provides an entry point to find what the misconception could be and how to address it. By doing this, we are promoting tam work. The material is something that concerns both the instructor and the student, and together they can work towards it. 

\subsubsection{Make it personal}

This pillar has to do with \emph{relevance}. A very common situation in Improv is when two performers spend their time talking about other people not present on stage, or about future or past actions. For the audience this is not relevant. They didn't go to see two people talking about stories or plans. It is not relevant for them to see two performers talking about other fictional characters. The audience mainly went to \emph{watch} a scene, not to listen a story. The audience is more interested in what the performers \emph{do}, what the dynamics among them are, what conflicts, relationships, or adventures the characters face in the plot. That is what is relevant for them. 

When performers start talking about people that are not present, or about actions in the past or the future, the scene can quickly turn very abstract and difficult to follow. Suddenly the audience is just listening to a story and not watching a performance. The issue then affects both the performers and the audience. For the performers, the risk is to become too abstract and for the audience, the risk is to disconnect from the scene. 

Likewise, while teaching we can fall into lack of relevance. When the contents are not relevant to either the instructor nor the students, engagement becomes a very difficult task. If the narrative and the context of the lesson are foreign to the students and to the instructor, no real motivation for the content can be obtained. To make it personal means to look for the relevance of the material for the students and for the instructor. 

The relevance might be practical, as direct concepts of applications that students will use later in their career. It can also be procedural, in developing necessary skills that will be relevant for other topics. This motivational component is important into building a reason why the students should be interested in learning the material and why the instructor is interested in facilitating it. 

A very successful way to make the material relevant is to connect it with applications important to the specific audience. When talking to a group of engineers, making it personal could mean to have examples of applications in engineering and not social sciences, for example. Even when questions are asked during class, a way to make it personal is to identify the interests of the student asking the question, and to direct the answer with these interests in mind. 

\subsubsection{Be present}

Being present is a pillar that has to do with \emph{consciousness}. During an Improv scene, performers need to be conscious about each other. They need to react to what the other performer does and says. In order to transmit a sense of unity, performers need to respond to the flow of the scene. If a performer does not react to something that happened in the scene, the audience will notice this and will sense disconnection on stage. Being present requires performers to react to what is being said or done, how it is said or done, and  when and where it happened. 

This focuses on details. When performers are so wrapped up into their own performance, it is very easy to overlook at the details happening on stage. This is a common mistake where performers pay more attention to themselves than to the scene. By paying attention to the present, performers are able to connect better with a common reality and to build a scene following the flow of the events that naturally happen on stage. 

For the instructor, this translates as awareness. By promoting awareness in the classroom, the instructor can more effectively generate and guide discussions. She would also be able to pace her lesson more effectively in order to accommodate students. A common mistake made by instructors is to teach only to themselves, ignoring and not reacting to what is happening in the classroom. By becoming more present, the instructor needs to be willing to deviate a bit from her original lesson plan in order to react to what might happen in class. This does not mean to change the entire lesson on the spot, but to incorporate the elements of class events into the discussion of her material. It also means to corroborate that she understood the students' questions. Many times as instructors, we interpret questions in one way, answering a different question than the original one.

\subsection{Rules of Improv}

The rules of Improv are guidelines to promote the pillars described in \Cref{sec:impro}. These practical techniques are helpful for giving  meta-structure to improvised scenes. Performers use rules to  approach the creative process in a more systematic way. Improvisation is a very free and open art form. Hence, it does not need to satisfy rules. The usage of the word \emph{rule} is somehow loose in this context. Perhaps a better expression would be \emph{Improv techniques}.

Traditioinally, Improv courses make the study of these rules the central part of their curriculum. Mastering these and having enough practice is the most common way to become an improviser. Usually, Improv is taught using games, exercises, and practicing short scenes that incorporate these rules. The particular approach is very different from place to place, but in general, this approach is universal. 

In the context of education, these rules are also referred as \emph{disciplined improvisation}\cite{sawyer2004creative, sawyer2011structure}. This means that we are controlling the flow of ideas while still keeping room for little deviations. 

Five of the most common Improv rules are: 
\begin{itemize}
\item Yes, and
\item Active listening
\item Embrace failure
\item Be present
\item Teamwork
\end{itemize}

These rules embody the pillars of Improv in a more practical way. They serve to ground the general Improv philosophy and make it more applicable. 

\subsubsection{Yes, and}

This is perhaps the most famous rule of Improv. ``Yes and'' means to acknowledge the other and to be present. By \emph{yes-anding}, performers are able to create a common reality and keep the flow of the scene. A tricky part of this rule of Improv is to distinguish between yes-and and becoming a yes-person. To yes-and means to acknowledge and agree with the reality and the reality of the others. It does not necessary mean to share a feeling, point of view, or goal. ``Yes, I recognize you are cold, but I am not,'' ``yes, I can see you are afraid of dogs, but I can be brave and face my fears,'' or `` yes, I acknowledge you want to go to Disney for the break, but I am not a big fan.'' In many cases, to yes-and also means to share the perspective, but the main focus of this rule is to recognize reality and not to deny it. 

When performers yes-and, they agree on a solid common context. This provides a good foundation to move the scene forward, without the need to fix plot holes or to spend extra time justifying events. This provides a more linear flow in the scene. When performers don't yes-and, it is very easy for them to fall into circles. The plot usually does not progress and the story feels stagnant.

For the instructor, specially during lesson, this rule becomes key. Yes-anding generates many teaching moments that can otherwise be lost. 

\begin{example}{Yes-and in class.}{1}
In one of my  precalculus lectures, we were discussing the idea of domain of a function. I gave the square function and asked the class to give examples about possible domains for this function. 

One student suggested that a possible domain was ``cosine.'' I was tempted to dismiss this answer, but went ahead and yes-anded their idea. I then asked the class if they agreed with that answer. Some students replied that this answer didn't work. I asked why and some said ``because it is not a number.'' Here, we were making the connection that domains had to do with numbers. So I wrote on the board, $$\cos(0)\,,$$ and asked again if that would be sufficient answer. Then some said that it was not enough because ``we needed a set of numbers.'' Someone suggested an interval, and then someone else suggested, $$\cos(0),\cos(1)\,.$$ Finally, I directed them to writing this using interval notation and following the increasing order of endpoints, for example
$$\left[\cos(1),\cos(0)\right)\,.$$
\end{example}

The situation described in Example \ref{ex::1} happened to me a couple of years ago in a large classroom with around 180 students. Even though I had a lesson plan for that day, I could not had prepared for a situation like this. Many times, as instructors, we simply dismiss answers that could seem out of place, and provide ourselves the answers to the questions. This would loose a valuable teaching moment. The student that replied with the \emph{apparent} out of context answer is telling the instructor that there is a misconception that needs to be addressed. Also, by yes-anding the answer, we show that we appreciate the input from students, that we don't expect perfect answers, and that we make them feel part of the learning experience. 

These interactions provide only a small deviation from the original lesson plan. The event outlined in Example \ref{ex::1} took place in the span of a couple of minutes. Also, it did not change the core lesson or content. In general, these interactions provide a way to connect with students and to take advantage of the teaching moments that naturally arise when interacting with them.

\subsubsection{Active listening}
Listening is more than just paying attention\cite{robertson2005active}. In Improv, listening plays a fundamental role since there is no script and the players have no notion of what to expect in a scene. Listening involves full and undivided attention\cite{robertson2005active}. By giving the other performer full attention, improvisers are able to fully grasp their understanding of a certain event, which in turns allows them to react and adjust accordingly. 

A performer needs to focus on what the other performer is saying, how she says it, when, what actions did she make, what reactions did the audience have, etc. When a performer is able to actively listen, she can become part of the story that is being created on stage. She gets integrated to the narrative as an active part of it. When this rule is not present, performers tend to become isolated from the reality created for the scene. This disconnection can produce a forced and unnatural flow of the plot in the scene. 

Something similar happens during a lesson. Communication can be understood as occurring when the behavior of one person affects the behavior of the other \cite{nvmandal}. As instructors, we seek to impact the behavior (learning) of our students. This includes verbal and non-verbal communication. Active listening includes both of these aspects \cite{doi:10.1111/j.1460-2466.2002.tb02559.x}. 

\begin{example}{Listening.}{2}
During one of my integral calculus lessons, we were finding the net change of the position of a particle with a given velocity. When setting up this example, I asked the students for ideas on how to solve the problem. One student said to ``find the equation between the initial and final points.'' This answer taken out of context might not reflect a good insight from the student. However, during the lesson we had solved a couple of similar problems where we \emph{found the integral between the initial and final points}. Also, the student was secure with his tone, and accompanied it with a hand motion describing an initial point and a final point with reference to the area under a curve.  In this case, the verbal plus the non-verbal cues, together with the lesson context gave me a good indication that the student had the right idea but \emph{said} the wrong thing.

I then addressed his idea and pointed out that a more accurate way to phrase it could be to ``find the integral between the initial and final points.''
\end{example}

By purposely focusing in all the communication components that the other person uses, we have more tools to better understand the standing of the other person, their mental and emotional context, which enables us to better respond and react to them \cite{johnson1999communication, doi:10.1177/009102609402300213}. By actively listening in Example \ref{ex::2}, the student could identify that his understanding was good but his terminology still needed to be refined. This is important for recognizing what specific misconception the student could have. That way, the efforts to clarify them can be better organized and targeted. During this situation, I could identify that the issue was the terminology and not conceptualization of the problem. When we as instructors don't actively listen to what happens in class, it is easy to mistakenly assume the student's misconception is something different than what it is in reality. 

\subsubsection{Embrace failure}

This rule is very important for Improv. Being too afraid of mistakes and failure can dramatically affect the performers creativity. Many beginners start by \emph{playing it safe}. They don't take risks, usually don't take the initiative in a scene, and over think their choices. This fearful approach is detrimental for Improv. One of the key things when performing is to be bold with one's choices. If a performer enters a scene saying she is cold and that she doesn't like cold weather, she need to keep this point of view. Even if the choices are more unexpected and unusual, she needs to keep her stance. If she started a scene rolling her Rs, then she needs to be bold and keep doing that. This could have happened out of a mistake. Maybe she mispronounced a word when starting. The audience doesn't need to know it was a mistake. By embracing this apparent failure, the performer is redefining the scene's reality and shaping her character. This mistake actually gave her a cue to build a character for the scene, \emph{a person that mispronounce every single word}. 

Mistakes are one of the best ways to produce comedy. A good amount of comedy is unintentional. It just happens. And the way to find it is  allowing for \emph{mistakes} to happen.

Similarly, in a classroom context, most teaching moments are unintentional. These often arise by \emph{mistakes}. Mistakes are a sign that learning is happening. In particular in mathematics, mistakes are often regarded as a big failure. The common belief is that in mathematics \emph{there is only one answer}, but we often forget to reinforce that there could be multiple ways to arrive to that answer. Not everyone thinks the same, and when a student pursues her own way of thinking, mistakes can happen. Making mistakes is also a sign that there is creative thinking involved. 

\newpage 
\begin{example}{Embrace failure.}{3}
In a lesson about the Net Change Theorem for my Integral Calculus course, I was working an example about a tank being filled up with water. The tank had capacity 10 gal, it was already filled up to 2 gal, and there is a source purring water into the tank at a rate of $$r(t)=e^{-t}\text{ gal/min}\,.$$ The question then was to find out how long did we have to wait for the tank to fill up.

Solving this problem requires to solve the equation, 

$$\int_0^T r(t)\,dt=10-2=8\,,$$
which leads to,

$$\left.-e^{-t}\right|_0^T=8\,,\quad T=-\ln(-7)\,.$$

This leads to complex valued solutions! While solving this example, I realized that I had made a mistake in the definition of the rate function. It originally was, $$r(t)=e^{t}\,.$$
I had made a typo when starting the example. When we reached to the last line of the example, I noted that this result would not make physical sense. Immediately after, a student asked if that meant that the tank would never be filled up, and this lead to a fruitful discussion. After this, I asked the students what changes could we make to the rate function so we actually be able to fill up the tank in finite time. Finally, I also noted that I had made a mistake on copying down the problem.
\end{example}

Even in scripted lessons, mistakes could happen. Instructors and students are human, and mistakes are expected. Situations like Example \ref{ex::3} can generate valuable teaching moments. During this, we were able to explore problems that didn't give a positive real solution. This gave us the chance to think and discuss what happens in those cases, and the effect of considering real life parameters, such as time. This also provided an opportunity to think about the general properties  a rate function should satisfy in order to lead to a valid solution. This is an indirect way to reinforce the conceptual understanding of the problem. Finally, embracing failure as an instructor also makes students to see mistakes as something normal in learning. Failure is a big stigma in our culture, but it is a fundamental part of learning. By normalizing and handling mistakes, we as instructors are setting the classroom culture that if someone makes a mistake, we can recognize that mistake and fix it. Even learn from it. This can take the burden from students for perfectionism, which can hinder their learning experience.

\subsubsection{Be present}

This is one of the most important features of Improv. When in stage, the audience can really connect with performers that develop a story that happens in the present time. Watching characters talk about the past or planning the future can make the audience bored and disengaged from the show. Both the performers and the audience are in the same physical location, therefore to achieve a better connection, they also should be in the same mental and emotional location.

Improv requires intense concentration and focus. When performers are not present, by thinking about life problems, their personal projects for the future, the bills they have to pay, or that they don't have food for breakfast the next day, they can miss important details that will make them fall out of the natural flow of the scene. 

Being present encloses many things which translate into the classroom. After being physically present, we need to be mentally present. One term comes to mind when talking about being mentally present: \emph{mindfulness}. The definition of mindfulness\cite{webster1999} states ``the practice of maintaining a nonjudgmental state of heightened or complete awareness of one's thoughts, emotions, or experiences on a moment-to-moment basis.'' Mindfulness has to do with the awareness of thoughts, emotions, and experiences. 

These three components are fundamental in the classroom. First, being mindful or aware of our thoughts. Here we have two aspects, as both the instructor and the students need to be mindful. This means, to be aware of one's state of mind. With the incorporation of more and more devices in our everyday lives, this process has become more challenging. In 2018, Americans looked at their phones on average 52 times a day\cite{deloitte2018}. Cellphone usage has proved to be a
detrimental factor in student performance \cite{doi:10.1177/2158244015573169, lepp2015exploring}. With all this distractions, becoming aware of our own thoughts is a challenge in the classroom. Being aware of our thoughts means to administer them according to the situations and circumstances we experience. Second, emotional awareness deals with the awareness of the other person. This entails knowing what the goals of students are, what is their motivation for being in class, and what is their current state in life. This awareness is important for connecting with them in an emotional way. This provides a better way for communicating the material in a more meaningful and effective way. Lastly, experiential awareness concerns with the reality of students and instructors. Do they work? Are they freshmen? Are they taking another math class? Motivation is an important component in the learning process \cite{doi:10.1080/00220670209596607, doi:10.3102/0013189X12459679}.


\begin{nonexample}{Be present.}{1}
While going through a proof for one of my Introduction to Number Theory courses, I found myself working the proof just by memory. I wrote the proof on the board and said the argument out loud, but was disconnected from it. 

My mind started drifting away from the proof and from the class. This produced a distance between the students and the lecture. Usually the class atmosphere was very interactive, and during this lecture it turned out dry and passive. 
\end{nonexample}

Situations like Non-Example \ref{nex::1} are not atypical for instructors. In certain situations it is easy to drift away from class. Sometimes as instructors we fall in the trap teaching class at ourselves, not at the students. This is the complete opposite of mindfulness in the classroom. 

Students are not immune to this lack of mindfulness. Sometimes life has its own issues and as humans, it is difficult to put those aside and to completely focus on class. As instructors, we cannot fully avoid students to drift away, but we can be strategic about first being present ourselves and second to promote students being present as well.

\subsubsection{Teamwork}

Another basic skill in Improv is teamwork. When performers create a scene, this is the result of everybody's work on stage. Everyone is part of the creative process and the execution of the scene. Even the audience has influence on how the scene turns out. 

Considering teamwork in an impro scene entails two things: the responsibility of supporting your scene partners and the certainty that they will look out for you. It is common to hear performers say to each other \emph{``got you back''} before going to stage. 

If a performer says on stage that her name is \emph{Ross} and another performer misheard and call her \emph{Rose}, a way to show teamwork is to the first performer to accept that her name from that point on is Rose. This provides trust, which makes performers focus on being creative and not being fearful of committing a mistake on stage. 

Something similar happens in the classroom. It is important to promote an atmosphere of teamwork in class. The instructor should strive to \emph{have her students' backs}. This can help diminish the fear of failure some students can have. 

\begin{example}{Teamwork.}{4}
During office hours for my Vector Calculus class, a student told me that she felt that some of the concepts in class where not completely clear for her and that homework had being very demanding. Nonetheless, she told me she didn't feel bad about it since the environment in class made her feel safe and confident that later she will understand the material, even if at that moment things where not completely clear. 
\end{example}

It is important for the instructor to convey the message that her role is of helping students. Even though this is assumed, in practice many students are not fully aware of this premise. Sometimes class feels like a battle between two teams, students and instructors, instead of both being part of the same team. The environment in Example \ref{ex::4} provides a set up where more progress can be made. When students and instructors feel part of the same team, both push towards the same direction and the effort required is less while the progression is more. If the perception is contrary, the efforts are mainly spend in the personal dynamics and not in the learning.

\section{Improv seminar}

Many current efforts to bring Improv skills and techniques involve the students directly\cite{berk2009whose, borko1989cognition, doi:10.1080/10511970.2012.754809}. These efforts are a great way to develop students' skills involved in active learning. Upfront, Improv games and strategies inside the classroom might seem disconnected from the mathematical content, but the emotional and motivational components developed become very important for active learning.

Improv skills can also be used in order to complement active learning with active teaching. Here, the focus is on the instructor's side. In this approach, Improv can be done outside of class by the instructors, leaving the students aside. This becomes even more relevant when the class size is big and where promoting Improv with students might be difficult to manage.

Pursing Improv practices with instructors or with the teaching team then becomes a Professional Development technique. This arises as a natural complement to active learning since it provides instructors with the necessary tools to manage and facilitate the active learning environment for students. 

One way we approached this professional development in the Mathematics Department at the University of California, Santa Cruz was through an Improv Seminar. This was aimed primarily at Teaching Assistants and Tutors. We ran weekly sessions were participants explored Improv concepts and discovered how these related to teaching. The Seminar was structured in an interactive way, were the session leader would introduce the concepts, facilitate discussion, and lead Improv exercises. 

Many Improv exercises are also good warm-up and ice-breaking techniques.  These are tools that can be used for class management for instructors, and for discussion or section management for teaching assistants and tutors. They also provide ways to practice public speaking and stage presence, which are fundamental for an instructor regardless of their teaching approach. 

\subsection{Improv Games}

Improv is usually classified as short or long form. Long form Improv seeks to create a scene with a plot while short form can be thought of as games. Games involve the key elements of Improv and usually focus in developing these as oppose to applying them to create a story. Many Improv courses start with a curriculum based on games as they are more approachable at first and usually are very engaging. 

The Seminar consisted of a total of eight meetings. We focused on different skills each meeting and explore different games that incorporate these skills. 

The aspects we covered were:
\begin{enumerate}[\bfseries{Session} 1:]
\item Yes, and
\item Be present
\item Teamwork
\item Failure
\item Association
\item Being obvious
\item Listening
\item Details
\end{enumerate}

The sessions and games were not completely disjoint. Some games involve more than one skill and provide a nice way of reinforcing certain techniques. 

\subsubsection{Yes, and}

As mentioned before, this skill involves \emph{accepting the other person}. Some selected games explored in the Seminar were:

\begin{itemize}
\item{\bf Yes and}: Participants go around a circle and each say a phrase one at a time. Each person must start by agreeing to the previous person's statement by saying ``yes, and,'' and stating something related to that previous statement. For example, if one persons says ``I love apples,'' the next person might say ``yes, and an apple a day keeps the doctor away.''
\item{\bf Tree}: In this game people form a circle. One person goes to the center and says ``I'm a tree.'' Another person then goes to the center and adds a related object. For example, ``and I'm a apple.'' Then a third person goes to the center and adds a final object to the scene, eg. ``and I'm Isaac Newton.'' Finally, the first person picks one of the two other objects and the chosen one stays to start a new scene.
\item {\bf Storytelling}: This game is similar to \emph{Yes, and} but now the goal is to create a story. That means that the new statements not only have to agree with the previous one, but also must be in line with all the other ones. The first person says ``once upon a time'' and adds a statement. Every person after that adds to the story.
\end{itemize}

\subsubsection{Be present}
\begin{itemize}
\item{\bf Woosh-bang-pow}: Players form a circle and pass around an invisible ball according to the following rules: pass the ball to your neighbor following the direction it came while saying out loud ``Woosh," or return the ball to your neighbor reversing the direction it came while saying ``Bang," or pass the ball skipping your immediate neighbor following the direction it came while saying ``Pow."
\item {\bf Enemy-defender}: Players start walking around the room in no particular direction. Secretly, every player chooses an \emph{enemy} that they need to hide from. Next, every player also secretly chooses a \emph{defender} that will protect them from their enemy. These roles are kept secret from everybody. When the game starts, every player should walk in such a way that their defender is between them and their enemy.
\item {\bf Last word}: In a circle, players go around saying phrases. Each new phrase must start with the last word the previous player used. This can also be done with words and letters.
\end{itemize}

\subsubsection{Teamwork}
\begin{itemize}
\item {\bf Mind meld}: Players form a circle. Two players say one word each and then at the same time, they try to say the same word, combining the two previous concepts. This goes around in a circle.
\item{\bf Count to 20}: With their eyes closed, players try to count up to 20 one number and one player at a time. If two players talk over each other, the count goes back to 0. 
\item {\bf The expert}: Three players sit in the front. They will form a single entity called \emph{the expert}. This entity is a world famous expert in a made up topic. A host will interview the expert by asking questions. The expert will then answer the questions. The three players must sequentially give one word at a time to build the answer. 
\end{itemize}

\subsubsection{Failure}
\begin{itemize}
\item{\bf Naming things}: Players wander around the room and start point at things. The first stage is to name the thing they are pointing. Do this for some time. The second stage is to name the last thing they pointed at while pointing a new one. Do this for some more time. The last stage is to name the second to last thing they pointed at. 
\item {\bf Five things}: In a circle, a player starts with the phrase ``name 5 things that" and completes it with any idea. Then their neighbor must come up with 5 things related to this idea, one at a time by saying ``One:" then the next idea ``Two:," etc. Then this goes in a circle. 
\item{\bf 3-7}: In a circle, players go around counting up from 1. When a player hits a multiple of 3, the player says ``Ping." When a player gets a multiple of 7 says ``Pong." When the number is a multiple of both, the player says ``Ping-Pong." The game is restarted when a player makes a mistake. 
\end{itemize}

\subsubsection{Association}
\begin{itemize}
\item{\bf Math mash up}: Players have to do a two minute presentation on a topic. The topic is chosen at random with the following strategy: there are three bins labeled ``Field," ``Problem," and ``Theorem." Each bin has a list of items. The player randomly selects one item per bin to form their topic. 
\item{\bf Association chain}: One player at a time is given two random words or concepts. The player then has to start from one words and give new terms that relate to the previous one until they end with the other word initially given. 
\item {\bf Translator}: Two player go to the front. One player is an alien and the other is the translator. A host will ask the alien question, and the translator with translate the questions and the answers. The translator should ask and reply incorporating the physical language used by the host and the alien. They also should interpret the gestures used and fit their translation accordingly. 
\end{itemize}

\subsubsection{Being obvious}
\begin{itemize}
\item {\bf Ding}: Two players go to the front and improvise a dialogue or a scene. When the host rings a bell, the players should repeat the last part they did but they should say a different thing. 
\item {\bf Barney}: In a circle, one player stand at the center and become "The caller." The caller says a letter and points to any player in the circle. That player must come up with a the name of a person, a service or object that can be sold, and a location, starting with the letter The Caller gave. If the player cannot come up with these, they become the new Caller. 
\item{\bf Presents}: In pairs, one player grabs an invisible present from the floor. Then the player starts describing the present and gives it to the other player. In return, this gives the first player another present. 
\end{itemize}

\subsubsection{Listening}
\begin{itemize}
\item{\bf Silent scenes}: Two players at a time hold a scene where they interact without words. 
\item {\bf Research Microscope}: One player starts explaining their research. The host will randomly say ``more" or ``less" depending if they want more or less details about the research.
\item{\bf Repeater}: Any group pf players do a scene. Every time one player wants to intervene, they need to repeat the entire dialogue the previous player did (without that player's  previous repetition). 
\end{itemize}

\subsubsection{Details}
\begin{itemize}
\item {\bf Small/Big}: One player starts a monologue in a random topic. The host will randomly say ``small" or ``big" depending if they want a smaller or bigger picture of the topic. 
\item{\bf Adjectives}: One player goes to the front an is given an object. The player starts describing this object as much as possible. The purpose is to describe it as long as possible.
\item {\bf Explain the math}: One player at a time will explain the math written on the board done by someone else. The purpose is to go into as many details as possible. 
\end{itemize}

\section{Conclusions}

Unexpected situations are actually expected in an active learning environment. Even the mpst experienced instructor will face unforeseen situations in her class for which she couldn't have prepare beforehand. Therefore, it is important to embrace these situations. These unexpected situations will generally arise from student-student or student-instructor interaction, generating natural teaching moments. During these moments, students are highly connected with the learning environment, providing a good avenue for teaching. 

A structured approach to manage these unexpected moments is provided by Improv, which despite not being scripted, supplies with a more abstract level of managing the unforeseen. Just as any other art form, Improv is something that needs training and practice. These same principles can be applied to the teaching context, in which the instructor can be trained to manage and utilize these unexpected situations to generate teaching moments. 

Improv practices and games can be used as a complement to active learning training. These skills become crucial for connecting, not only on an intellectual, but also on an emotional level with students. 

\bibliographystyle{plain}
\bibliography{references}

\end{document}